\documentclass[makeidx]{article}
\usepackage[latin1]{inputenc}
\title{Upper bounds for the ground state energy of the Laplacian with zero magnetic field on planar domains}
\author{Bruno Colbois and Alessandro Savo}
\date{\today}
\pdfoutput=1
\usepackage{makeidx}
\usepackage{amssymb , amsmath, amsthm,graphicx, float}
\usepackage{graphicx}
\usepackage{xcolor}
\newtheorem{defi}{Definition}
\newtheorem{theorem}[defi]{Theorem}

\newtheorem{example}[defi]{Example}

\newtheorem{lemme}[defi]{Lemma}
\newtheorem{cor}[defi]{Corollary}

\newcommand{\R}{\mathbb R}

 \newcommand{\twosystem}[2]{\left\{\begin{aligned} &#1\\ &#2\end{aligned}\right.}
\newcommand{\threesystem}[3]{\left\{ \begin{aligned}&#1\\ &#2\\&#3\end{aligned}\right.}

\newcommand{\nero}{\smallskip$\bullet\quad$\rm}

\newcommand{\abs}[1]{\lvert{#1}\rvert}

\newcommand{\reals}{{\bf R}}

\newcommand{\real}[1]{{\bf R}^{#1}}

\newcommand{\bd}{\partial}

\newcommand{\matrice}{\begin{pmatrix}}
\newcommand{\ok}{\end{pmatrix}}

\begin{document}
\maketitle

\maketitle
\abstract{We obtain upper bounds for the first eigenvalue of the magnetic Laplacian  associated to a closed potential $1$-form (hence, with zero magnetic field) acting on complex functions of a planar domain $\Omega$, with magnetic Neumann boundary conditions. It is well-known that the first eigenvalue  is positive whenever the potential admits at least one non-integral flux.  By gauge invariance the lowest eigenvalue is simply zero if the domain is simply connected; then, we obtain an upper bound of the ground state energy depending only on the ratio between the number of holes and the area; modulo a numerical constant the upper bound is sharp and we show that in fact equality is attained (modulo a constant) for Aharonov-Bohm-type operators acting on domains punctured at a maximal $\epsilon$-net. In the last part we show that the upper bound can be refined, provided that one can transform the given domain in a simply connected one by performing a number of cuts with sufficiently small total length; we thus obtain an upper bound of the lowest eigenvalue by the ratio between the number of holes and the area, multiplied by a Cheeger-type constant, which tends  to zero when the domain is metrically close to a simply connected one. 

\bigskip

\noindent {\it Classification AMS $2000$}: 58J50, 35P15\newline
{\it Acknowledgments:} Research partially supported by INDAM and GNSAGA of Italy

%\tableofcontents
\large
 
\noindent

\section{Introduction}
We consider a smooth, connected and bounded domain $\Omega \subset \R^2$ of area $\vert \Omega \vert$. Let $A$ be a closed $1$-form and introduce the magnetic Neumann Laplacian $\Delta_A$ with potential $A$ acting on functions $u \in C^{\infty}(\Omega,\mathbb C)$. It is the operator $\Delta_A=(\nabla^A)^{\star}\nabla^A$ where the connection $\nabla^A$ is defined as $\nabla^Au=\nabla u-iuA^{\sharp}$ and  $A^{\sharp}$ is the vector potential, dual to the $1$-form $A$.
The following notation is sometimes used:
$$
\Delta_A=(i\nabla+A^{\sharp})^2,
$$

We take  (magnetic)  Neumann boundary conditions and then we study the eigenvalue problem:
\begin{equation}\label{eproblem}
\twosystem
{\Delta_Au=\lambda u\quad\text{on}\quad \Omega}
{\nabla^A_Nu=0\quad\text{on}\quad \bd\Omega}
\end{equation}
where $N$ is the inner unit normal. We are interested in the ground state energy (lowest eigenvalue) $\lambda_1(\Omega,A)$, which is the minimum of the Rayleigh quotient:
\begin{equation}\label{minmax}
\lambda_1(\Omega,A)=\min\Big\{\dfrac{\int_{\Omega}\abs{\nabla^Au}^2}{\int_{\Omega}\abs{u}^2}: u\in H^1(\Omega)\setminus\{0\}\Big\}.
\end{equation}

If $A=0$ the spectrum of $\Delta_A$ coincides with the spectrum of the usual Laplacian under Neumann boundary conditions. The same is true when $A$ is an exact one-form, thanks to  the well-known {\it gauge invariance} of the magnetic Laplacian, which implies in particular that:
\begin{equation}\label{gi}
\lambda_1(\Omega,A)=\lambda_1(\Omega,A+df)
\end{equation}
for any smooth function $f$. The two-form $B\doteq dA$ is called the {\it magnetic field} associated to the potential $A$.
It turns out that, even if the magnetic field is $0$, the (closed) potential $A$ can affect the ground state energy: this is related to a phenomenon in quantum mechanics called {\it Aharonov-Bohm effect}. To see this, we introduce the \emph{flux} of $A$ across the closed curve (loop)  $c$ as the quantity:
$$
\Phi^A_c=\frac{1}{2\pi}\oint_cA
$$
 (we don't specify the orientation of the loop, as it will be irrelevant for our bounds).

\medskip

It turns out that $\lambda_1(\Omega,A)=0$ if and only if $A$ is closed and the cohomology class of $A$ is an integer, that is, the flux of $A$ around any loop is an integer. 
 
\smallskip
	 
This fact was first observed by Shigekawa \cite{Sh} for closed manifolds, and then proved in Theorem 1.1 of  \cite{HHHO} for manifolds with boundary. This remarkable feature of the magnetic Neumann Laplacian shows its deep relation with the topology of the underlying domain $\Omega$. 
For a more detailed introduction to the magnetic Neumann Laplacian associated to a closed potential, see the introduction of \cite{CS2} and the references therein.

\smallskip

It is precisely the goal of this note to investigate how the topology and the geometry of the domain $\Omega$ influence the ground state energy $\lambda_1(\Omega,A)$ when the magnetic field is zero. Therefore, from now on, unless otherwise stated:

\nero  the potential $A$ in this paper will always be a {\it closed} one-form. 

\smallskip

Note that, in view of Shigekawa's remark, any lower bound of the  ground state energy should somewhat depend on the distance of the fluxes of $A$ to the lattice of integers which, for a single loop $c$, is defined as:
$$
d(\Phi^A_c,{\bf Z})\doteq\min\{\abs{\Phi^A_c-k},\, k\in{\bf Z}\}.
$$
In our previous papers \cite{CS1} and \cite{CS2} we obtained lower bounds for the ground state energy. 
In \cite{CS1}  we proved  an estimate from below of the first eigenvalue of a Riemannian cylinder; applied to a plane annulus $\Omega=F\setminus \bar G$, with $F$ and $G$ smooth and convex, the lower bound becomes:
\begin{equation}\label{lbannulus}
\lambda_1(\Omega,A)\geq\dfrac{4\pi^2}{\abs{\bd F}^2}\dfrac{\beta^2}{B^2}d(\Phi,{\bf Z})^2
\end{equation}
where $\Phi$ is the flux of $A$ across the inner boundary component $\bd G$ and $\beta$ (resp. $B$) is the minimal length (resp. the maximal length) of a segment contained in $\Omega$ and hitting the inner boundary $\bd G$ orthogonally. We  call $\beta$ and $B$ the {\it minimal}, resp. {\it maximal} width of $\Omega$, respectively; obviously $\beta$ is also the minimum distance between the inner and outer boundary curves. 

In the subsequent paper \cite{CS2} we improved the result to get a lower bound depending on $\frac{\beta}{B}$, rather than $\frac{\beta^2}{B^2}$ (see Theorem 1):
\begin{equation}\label{sharp}
\lambda_1(\Omega,A)\geq \dfrac{\pi^2}{8}\dfrac{\abs{F}^2}{\abs{\bd F}^2D(F)^4} \frac{\beta(\Omega)}{B(\Omega)}\cdot d(\Phi^A,{\bf Z})^2.
\end{equation}
where $D(F)$ is the diameter of $F$; this linear dependence is in fact sharp, as shown in \cite{CS2}. We will in fact use \eqref{sharp} in Section \ref{example}, formula \eqref{minoration}.
In \cite{CS2} we also extend the lower bound to domains with an arbitrary number of holes. 

\smallskip

Upper bounds for the spectrum of the magnetic Schr\"odinger operator, for an arbitrary potential one-form $A$, where considered in \cite{CEIS}. Some of them are consequence of the inequality
$
\lambda_1(\Omega,A)\leq \mu_1(\Omega),
$
where $\mu_1(\Omega)$ is the lowest eigenvalue of the Schr\"odinger operator $\Delta+\abs{A}^2$, with Neumann boundary conditions. In particular, Theorem 3 in \cite{CEIS} gives an upper bound of the first eigenvalue when the potential is a harmonic one-form, which depends on the volume of $\Omega$ and the distance (taken in $L^2$) of $A$ to the lattice of  integral harmonic one-forms. However, this upper bound is difficult to compute, in general.

\smallskip

The scope of this paper is to prove upper bounds of $\lambda_1(\Omega,A)$ which are computable, and depend explicitly on the topology and the geometry of $\Omega$. 
The topology of a planar domain $\Omega$ is specified by the number $n=n(\Omega)$ of holes, and in fact our first main result, Theorem \ref{flux}, gives an upper bound of the ground state energy depending only on the area of $\Omega$ and the number of holes; up to a numerical constant, the bound is sharp, and is achieved for a certain class of punctured domains (see Theorem \ref{epsilonnet}).   Note that, if $n=0$, 
$\Omega$ is simply connected, $A$ is exact and then $\lambda_1(\Omega,A)=0$: one could  intuitively argue that if it is possible to transform a domain $\Omega$ into a simply connected domain by deleting a family of segments of small total length then $\lambda_1(\Omega,A)$ should be small. We somewhat show that in fact this is the case in our second main result, Theorem \ref{main2}.

We now give the precise statements of our results.

\section{Main results} We have already remarked that $\lambda_k(\Omega,0)$ is just the $k$-th eigenvalue of the usual Laplacian with Neumann boundary conditions. Note that $\lambda_1(\Omega,0)=0$, the associated eigenspace being one-dimensional, spanned by the constant functions. Then, one could ask if $\lambda_1(\Omega,A)$ could be somewhat compared to the first positive Neumann eigenvalue, that is, to $\lambda_2(\Omega,0)$ (but in fact we will see that there is no a priori inequality between the two eigenvalues, see below).  To that end,  recall the Szeg\"{o}-Weinberger inequality, stating that the Neumann ground state is bounded above by that of the disk $\bar\Omega$ of the same area:
$$
\lambda_2(\Omega,0)\leq \lambda_2(\bar\Omega,0),
$$
which leads to an upper bound only in terms of the area:
\begin{equation}\label{weaksw}
\lambda_2(\Omega,0)\leq \dfrac{C}{\abs{\Omega}}
\end{equation}
where $C=\pi\lambda_2(B,0)$ and $B$ is the unit ball in $\real 2$.

Our first question was to see if a weak Szeg\"{o}-Weinberger inequality could possibly hold for $\lambda_1(\Omega,A)$ in this context: that is, can we find an absolute constant $C$ such that for every closed potential $A$ on $\Omega$ one has:
\begin{equation}\label{mw}
\lambda_1(\Omega,A)\leq \dfrac{C}{\abs{\Omega}}?
\end{equation}
A bit surprisingly, we find out that \eqref{mw} cannot hold in that generality; the estimate must in fact depend on the {\it topological complexity} of $\Omega$, that is, on the number $n=n(\Omega)$ of holes.

\begin{theorem}  \label{flux} Let $\Omega \subset \mathbb R^2$ be a bounded domain with smooth boundary having  $n$ holes. Then, there exists a universal constant $C>0$ such that for every closed potential $A$ we have:

\begin{equation} \label{borneflux}
\lambda_{1}(\Omega, A) \le C \frac{n}{\vert \Omega\vert}.
\end{equation}

One could take $C=544\pi$.
\end{theorem}

The constant $C$ is not optimal, but modulo a universal constant  the upper bound is sharp, as the next result will show. In other words there are examples of pairs $(\Omega,A)$ with $\Omega$ of fixed area, whose first eigenvalue grows proportionally to the number of holes $n$.

\subsection{Punctured domains and maximal $\epsilon$-nets} We introduce {\it punctured domains}: these are obtained by deleting  $n$ given points $\mathcal P=\{p_1,\dots,p_n\}$ from a given domain $\Omega$. We define:
$$
\lambda_1(\Omega\setminus\mathcal P,A)=\liminf_{\eta\to 0}\lambda_1(\Omega\setminus\mathcal P(\eta),A)
$$
where $\mathcal P(\eta)$ is the $\eta$-neighborhood of $\mathcal P$ (it obviously consists of a finite set of closed disks of radius $\eta$). It is not our scope in this paper to investigate the convergence in terms of $\eta$.

\medskip

A general lower bound for the first eigenvalue of punctured domains is given in Theorem 3 of \cite{CS2}. An interesting feature of punctured domains, which we will explicit in Section \ref{partition}  and which does not follow trivially from \cite{CS2},  is that their first eigenvalue could grow proportionally to the number of punctures, provided that the configuration $\mathcal P$ is a {\it maximal $\epsilon$-net}, which we are going to define.

\begin{defi} Given a convex domain $\Omega \subset \R^2$ with smooth boundary and  a number $\epsilon>0$, a maximal collection of points $\mathcal P_{\epsilon}=\{p_1,\dots,p_n\}$ with the following properties:

\nero $d(p_j,p_k)\geq \epsilon$ for all $j\ne k$,

\nero $d(p_j,\bd\Omega)\geq \epsilon$ for all $j$ 

\smallskip

\noindent {\it is called a {\rm maximal $\epsilon$-net}.}

\end{defi}

One should think of a maximal $\epsilon$-net as an optimal way of distributing a set of points inside $\Omega$ with the constraint of being at distance at least $\epsilon$ among themselves and at distance at least $\epsilon$ to the boundary.
Consider the harmonic $1$-form $A$ on $\Omega\setminus \mathcal P_{\epsilon}$ which has the same flux $\Phi>0$ around each of the holes $p_1,\dots,p_n$.  We denote by
$$
d(\Phi,{\bf Z})=\min\{\abs{\Phi-k}: k\in\bf Z\}
$$
the distance of the common flux $\Phi$ to the lattice of integers.
We then have:

\begin{theorem} \label{epsilonnet} If $\bd\Omega$ satisfies the $\delta$-interior ball condition, then, for all $\epsilon<\delta$ and for all maximal $\epsilon$-nets $\mathcal P_{\epsilon}$, one has:
$$
\lambda_1(\Omega\setminus \mathcal P_{\epsilon},A)\geq\dfrac{1}{64} \dfrac{{d(\Phi,{\bf Z}})^2}{\epsilon^2}
$$
 
In terms of the number of points $n=n(\epsilon)$ (hence, the number of holes), we have:
\begin{equation}\label{wein}
\lambda_1(\Omega\setminus \mathcal P_{\epsilon},A)\geq \dfrac{\pi}{256}\cdot\dfrac{n}{\abs{\Omega}}\cdot d(\Phi,{\bf Z})^2.
\end{equation}
\end{theorem}

%Let $\mathcal P_{\epsilon}=\{p_1,\dots,p_n\}$ be a maximal $\epsilon$-net. 

The strategy of the proof is to partition the given punctured domain in a family of convex domains with only one puncture, and then to apply a lower bound proved in \cite{CS1} to each piece of the partition. 

Recall that $\Omega$ satisfies the \emph{$\delta$-interior ball condition} if, for any $x\in\bd \Omega$, there exists a ball of radius $\delta$ tangent to $\bd \Omega$ at $x$ and entirely contained in $\Omega$. This is equivalent to saying that the injectivity radius of the normal exponential map is at least $\delta$; hence any point of a segment hitting the boundary orthogonally at $p\in\bd\Omega$ minimizes the distance to the boundary up to distance $\delta$ to $p$.

More generally, the inequality holds with ${\cal F}^2$ replacing $d(\Phi,{\bf Z})^2$, where
$$
{\cal F}^2=\min_{j=1,\dots, n}d(\Phi_j,{\bf Z})^2
$$
and $\Phi_j$ is the flux of $A$ around $p_j$.

\smallskip

Assuming constant flux $\frac 12$ around every point of the net, we see that the domain $\Omega\setminus \mathcal P_{\epsilon}$ (having $n$ holes and area $\abs{\Omega}$) satisfies the bounds:
$$
\dfrac{1}{1024}\leq \lambda_1(\Omega\setminus\mathcal P_{\epsilon},A)\cdot\dfrac{\abs{\Omega\setminus\mathcal P_{\epsilon}}}{n}\leq 544,
$$
showing that 
\eqref{borneflux}, modulo a  numerical constant, is sharp. 

\smallskip

A final question in this regard is the following:

\nero{\it  Is there an inequality relating  $\lambda_1(\Omega,A)$ with $\lambda_2(\Omega,0)$?}

\smallskip

The answer is negative. To show this, first consider that when $\epsilon$ is sufficiently small and $0<\eta<\frac{\epsilon}{2}$, one has, in the previous notation:
$$
\lambda_1(\Omega\setminus\mathcal P_{\epsilon}(\eta),A)>\lambda_2(\Omega\setminus\mathcal P_{\epsilon}(\eta),0)
$$
where $\mathcal P_{\epsilon}(\eta)$ is the $\eta$-neighborhood of $\mathcal P_{\epsilon}$. In fact, as $\epsilon\to 0$, the left-hand side diverges to infinity while the right-hand side  is uniformly bounded above by the Szeg\"{o}-Weinberger inequality
\eqref{weaksw}. 

In the other direction, remove from a fixed rectangle $F$ in the plane another smaller rectangle $G$ with fixed sides parallel to those of $F$, such that the boundary components of $F$ and $G$ get $\epsilon$-close to each other: see Figure 2  in \cite{CS2}. There it is proved that the resulting domain $\Omega_{\epsilon}$ is such that 
$
\lambda_1(\Omega_{\epsilon}, A)
$
converges to zero proportionally to $\epsilon$, where $A$ is the closed potential having flux $\frac 12$ around the inner curve $\bd G$.  Nevertheless, one observes that the Cheeger constant of $\Omega_{\epsilon}$ is uniformly bounded below by a positive constant $C$, which implies that $\lambda_2(\Omega_{\epsilon},0)\geq C>0$. Therefore, for $\epsilon$ small one has $\lambda_1(\Omega_{\epsilon},A)<\lambda_2(\Omega,0)$.

\smallskip

We remark that in \cite{FH2} the authors investigate the validity of the 
Szeg\"{o}-Weinberger inequality when the magnetic potential is non-zero (in particular, has constant norm). 
%%%

\subsection{An upper bound by a Cheeger type constant} First observe that, if $\Omega$ has $n$ holes, one  can suitably delete $n$ segments from $\Omega$ (joining different connected components of the boundary) and get a simply connected domain. We will establish an upper bound of $\lambda_1$ depending on the sum of the lengths of these segments, denoted by $h(\Omega)$. On one side, it will show that if $h(\Omega)$ is small enough, then the upper bound we get is better then the bound of Theorem \ref{flux};  on the other hand, we will construct an example showing that even if $h(\Omega)$ goes to $0$, $\lambda_1 (\Omega,A)\vert \Omega\vert$ may be large (and therefore the number of holes must be large).

\begin{defi} \label{Cheeger}Let $\Omega \subset \mathbb R^2$ be a bounded domain. An {\rm admissible cut} of $\Omega$ is a collection of segments $\Gamma=\{\Gamma_1,\dots,\Gamma_n\}$ such that $\Omega\setminus\Gamma$ is simply connected. Introduce the constant $h(\Omega)$:
$$
h(\Omega)= \min \{\sum_{i=1}^n h_i: \ h_i={\rm length}(\Gamma_i)\}
$$ 
where $\Gamma$ is a admissible cut of $\Omega$.
\end{defi}

The constant $h$ may be seen as an adapted Cheeger constant to measure how the topology (the number $n=n(\Omega)$ of holes) and the geometry (the lengths $h_j$ of the segments $\Gamma_j$) interact in order to affect the first eigenvalue $\lambda_1(\Omega,A)$. A natural question is for example to ask how small $h$ must be in comparison with $n(\Omega)$ in order to guarantee that $\lambda_1(\Omega,A)$ is uniformly bounded for a family $\Omega$ of domains of given area with a fixed number of holes. 

%This is the reason why we will add some technical constrains about $h$ and $h_i$ in the statement of Theorem \ref{main2} which are convenient for the proof and are satisfied when $h$ is small enough (the significant case).

\begin{theorem} \label{main2}
 Assume  that  $h(\Omega) \le   \frac{\vert \Omega\vert}{2\pi} $ and $h_j\leq 1$ for each $j=1,...n$. Then:
\begin{equation} \label{in4}
\lambda_1(\Omega,A)\le \frac{8\pi n(\Omega)}{\abs{\Omega}} 
\sum_{j=1}^n \frac{1}{\abs{\ln \frac{h_j}{2}}}
\end{equation}
where $h_j$ denotes the length of the j-th segment $\Gamma_j$ associated to $h$ and $n(\Omega)$ is the number of holes.
\end{theorem}

Note that we assume bounds on $\frac{h}{\abs{\Omega}}$ and on every $h_j$: this is a technical fact needed for the proof. On the other hand, Theorem \ref{main2} is meaningful and improves Theorem \ref{flux} in the special situation  where $h_j$ and $h$ are very small; the general situation is treated in Theorem \ref{flux}, which does not follow from Theorem \ref{main2}.

\begin{cor}
In particular, assume that $\Omega$ is doubly-connected. If 
 $h \le \min\{1,\frac{\abs{\Omega}}{2\pi}\}$, then
\begin{equation} \label{in2}
\lambda_1(\Omega,A) \le \frac{8\pi}{\vert \Omega \vert \abs{\ln\frac{h}{2}}}
\end{equation}
\end{cor}  

Note that for doubly connected domains one has $h=\beta$ where 
$\beta$ is the minimal width of $\Omega$, and also the minimum distance between the two boundary components. The corollary shows that if $\abs{\Omega}$ is fixed and the boundary components get very close (that is, $h=\beta$ tends to zero) then $\lambda_1$ tends to zero, which indeed improves Theorem \ref{flux}.   An interesting question is to see if the rate at which this happens, that is $1/\abs{\ln \frac h2}$, is actually sharp or can be improved. 

\medskip

When there is more than one hole, it is still possible to have an upper bound directly in terms of $h(\Omega)$. 

\begin{cor} \label{corh} Assume in addition that, in the definition of $h$, every $h_j\leq e^{-2}$. Then we have:
\begin{equation} \label{in3}
\lambda_1(\Omega,A)\le \frac{8\pi n(\Omega)^2}{\abs{\Omega}}\dfrac{1}{\abs{\ln(\frac{h(\Omega))}{n(\Omega)})}}.
\end{equation}
\end{cor}

For example, if $\Omega$ has area $1$ and $n$ holes, in order to guarantee that $\lambda_1(\Omega,A) \le 1$, we need to impose
$h(\Omega) \le n e^{-8\pi n^2}$.
 
\medskip
It is natural to ask what occurs when $h \to 0$ for domains of given area.
Clearly, if $n(\Omega)$ is fixed and $h\to 0$, inequality (\ref{in3}) implies that $\lambda_1(\Omega,A) \to 0$. 

\medskip
However, if $n(\Omega)$ is not fixed, the assumption $h\to 0$  does not imply that the first eigenvalue tends to zero. In fact, we can have $h$ arbitrarily small and, at the same time, $\lambda_1(\Omega,A) \vert \Omega \vert$ as large as one wishes. 
The next example is an illustration of this fact.

\begin{example} \label{example2} There exists a family of domains $\{\Omega_k\}_{k\geq 1}$ with area $\vert \Omega_k\vert \ge 1$ and with a fixed potential $A$ of equal flux $\Phi>0$ around each hole such that
$$
h(\Omega_k)\leq \dfrac{2}{\sqrt k}
$$
and, at the same time: 
$$
\lambda_1(\Omega_k,A)\geq c\sqrt k d(\Phi,{\bf Z})^2
$$

with $c=\frac{\pi^2}{2^{15} \sqrt{2}}$. The number of holes of $\Omega_k$ is $n(\Omega_k)=k^2$ (hence, it grows with $k$).
\end{example}

\section{Proof of Theorem \ref{flux}} Recall that we want to show that if $\Omega \subset \mathbb R^2$ is a bounded domain with smooth boundary having $n$ holes then, for every closed potential $A$, we have:
\begin{equation}  
\lambda_{1}(\Omega, A) \le  544\pi\frac{n}{\vert \Omega \vert}.
\end{equation}

\bigskip
\noindent
\textbf{Proof.}
The proof consists in three steps. First, using gauge invariance we replace the given potential $A$ by a new potential having the same flux but with poles at a certain collection of points $\{p_1,....,p_n\}$. The two corresponding magnetic Laplacians are unitarily equivalent and have the same spectrum. In the second step, we show the existence of a ball $B(p,r)$ of radius 
$r=\dfrac{1}{4\sqrt{\pi}}\Big(\dfrac{\abs{\Omega}}{n}\Big)^{\frac 12}$
 such that for each $i=1,...,n$, $p_i \not \in B(p,2r)$. Moreover we get a control of the area growth by the relation 
$$
\dfrac{\abs{B(p,2r)\cap\Omega}}{\abs{B(p,r)\cap\Omega}}\leq 34.
$$

\noindent In the last step, the fact that $p_i \not \in B(p,2r)$ for any $i$ will imply that $A$ is exact on $B(p,2r)$, hence, thanks to the control of the volume growth of the balls, one can control $\lambda_1(\Omega,A)$ by a standard cut-off argument for the usual Laplacian.

\medskip
\noindent
\textbf{Step 1.} The domain $\Omega$ is bounded by an outer closed curve $\Sigma_0$ and $n$ closed inner curves $\Sigma_1,...,\Sigma_n$. We assume that our closed potential $A$ has flux $\Phi^A_{i}$ around $\Sigma_i$. 

 \medskip
 \noindent
We choose $n$ points $p_1,...,p_n$ so that $p_i$ is inside $\Sigma_i$, and we write $(a_i,b_i)$ for the coordinates of $p_i$. Let $A_i$ be the $1$-form
$$
A_i(x,y)=\Phi^A_{i}\left( \frac{-(y-b_i)}{(x-a_i)^2+(y-b_i)^2}dx+\frac{(x-a_i)}{(x-a_i)^2+(y-b_i)^2}  dy\right).
$$
The flux of $A_i$ is equal to $\Phi^A_{i}$ around $\Sigma_i$ and it is $0$ around $\Sigma_j$ for $j \not =i$ (we assume that every $\Sigma_i$ is travelled once).
This implies that the fluxes associated to the $1$-form $\tilde A\doteq A_1+...+A_n$ are equal to the fluxes of $A$, and therefore $A-\tilde A$ is exact. By Gauge invariance, the operators $\Delta_A$ and $\Delta_{\tilde A}$ have the same spectrum, so that it suffices to find an upper bound for 
$\lambda_1(\Omega, \tilde A)$.

\medskip
\noindent
\textbf{Step 2.} First we prove the following estimate that will be used in the proof of Lemma \ref{ball} below. Let $a>b$ and consider the maximal number $N=N(a,b)$ of points at distance at least $b$ from one another which are in a ball $B(a)$ of radius $a$. We find:
\begin{equation}\label{estimateofn}
\Big(\dfrac{a}{b}\Big)^2\leq N\leq \Big(\dfrac{2a+b}{b}\Big)^2.
\end{equation}
To see this estimate, we denote by $x_1,...,x_N$ a maximal net of points at mutual distance  at least $b$ in the ball $B(a)$. The balls of center $x_i$ and radius $\frac{b}{2}$ are disjoint and contained in the ball of radius $a+\frac{b}{2}$, so that
$$
\sum_{i=1}^N\vert B(x_i,\frac{b}{2})\vert \le \vert B(a+\frac{b}{2})\vert
$$
which means that
$
N(a,b)\frac{b^2}{4} \pi \le \pi\Big(\frac{2a+b}{2}\Big)^2,
$
and then:
$$
N(a,b)\le \Big(\frac{2a+b}{b}\Big)^2.
$$
On the other hand, by maximality, the union of the balls $B(x_i,b)$ covers $B(a)$, and
$$
\vert B(a) \vert \le \sum_{i=1}^N\vert B(x_i,b)\vert
$$
so that 
$
\pi a^2\le N(a,b)\pi b^2
$
and $N(a,b) \ge \frac{a^2}{b^2}$ as asserted. This proves \eqref{estimateofn}.

\medskip
In order to construct an upper bound for $\lambda_1(\Omega,\tilde A)$, we will construct a test function with Rayleigh quotient $\le C\frac{n}{\vert \Omega\vert}$. This test function will be constructed geometrically with ideas coming from \cite{CM}, but much easier to apply in our case, because we are concerned only with the first eigenvalue.
Fix the number:
$$
r\doteq\dfrac{1}{4\sqrt{\pi}}\Big(\dfrac{\abs{\Omega}}{n}\Big)^{\frac 12}.
$$
Then we have the following fact.

\begin{lemme} \label{ball}There exists a point $p\in\Omega$ such that $p_j \not \in B(p,2r)$ for every $j=1,...,n$, and moreover
$$
\dfrac{\abs{B(p,2r)\cap\Omega}}{\abs{B(p,r)\cap\Omega}}\leq 34.
$$
\end{lemme} 

\noindent {\bf Proof of the lemma.}  
\medskip
From the definition of $r$ one sees that
$$
\sum_{j=1}^n\abs{B(p_j,2r)}\leq n\pi (2r)^2\leq 4\pi n\dfrac{1}{16\pi}\abs{\Omega}\dfrac{1}{n}=\dfrac{\abs{\Omega}}{4}.
$$
Set
$$
\Omega_0=\Omega\setminus \cup_j B(p_j,2r),
$$
so that
\begin{equation}\label{zero}
\abs{\Omega_0}\geq \dfrac 34\abs{\Omega}.
\end{equation}
For any $q\in \Omega_0$ we have clearly $p_j\notin B(q,2r)$ for all $j$. Take a maximal $r$-net in $\Omega_0$, say $\mathcal N=\{q_1,\dots,q_m\}$, so that $d(q_i,q_j)\geq r$ for all $i,j$ and by the  maximality of the net:
$$
\Omega_0\subset \cup_{j=0}^mB(q_j,r).
$$
This implies that
\begin{equation}\label{one}
 \sum_{j=1}^m\abs{B(q_j,r)\cap\Omega)} \ge \abs{\Omega_0} \ge \frac{3}{4} \abs{\Omega}.
\end{equation}

By the estimate \eqref{estimateofn}, for any $q\in\Omega$ the cardinality of the set $\mathcal N\cap B(q,2r)$ is at most $(\frac{2a+b}{b})^2$, with $a=2r$ and $b=r$, that is:
$$
\abs{\mathcal N\cap B(q,2r)}\leq 25.
$$
In other words, every point $q\in\Omega$ is in at most $25$ balls of radius $ 2r$ 
centered at a point of the net, hence
\begin{equation}\label{two}
\sum_{j=0}^m\abs{B(q_j,2r)\cap\Omega}\leq 25\abs{\Omega}.
\end{equation}
We can now prove that there exist $q_j$ such that
$$
\dfrac{\abs{B(q_j,2r)\cap\Omega}}{\abs{B(q_j,r)\cap\Omega}}\leq 34.
$$
Assume not. Then:
$$
\abs{B(q_j,2r)\cap\Omega}>34 \abs{B(q_j,r)\cap\Omega}
$$
for all $j$. We would then have, by \eqref{zero}, \eqref{one} and \eqref{two}:
$$
\begin{aligned}
25\abs{\Omega}&\geq\sum_{j=0}^m\abs{B(q_j,2r)\cap\Omega}\\
 &> 34 \sum_{j=1}^m\abs{B(q_j,r)\cap\Omega}\\
&\geq 34 \cdot  \frac{3}{4} \abs{\Omega}\\
& > 25\abs{\Omega}.
\end{aligned}
$$
which is a contradiction. The lemma is then proved. 

\medskip

{\bf Step 3.} We take a ball $B(p,2r)$ as in Lemma \ref{ball}. Then we can conclude as follows.  First, the restriction of $\tilde A=A_1+...+A_n$ to $B(p,2r)$ is exact, because the poles $p_1,...,p_n$ are not contained in the ball. Up to a Gauge transformation, we can replace the magnetic Laplacian  $\Delta_{\tilde A}$ by the usual Laplacian on $B(p,2r)$. 

We define a  function $u:\Omega\to\reals$ as follows:
$$
u(x)=\twosystem
{1\quad\text{if}\quad d(p,x)\leq r}
{-\frac 1rd(p,x)+2\quad\text{if}\quad d(p,x)\geq r}
$$
Note that $u$ is indeed supported on $B(p,2r)$; extending it to zero on the complement of the ball, we get a well-defined test function. As $\abs{\nabla u}\leq \frac 1r$, we see:
$$
\int_{\Omega}\abs{\nabla u}^2\leq \dfrac 1{r^2}\abs{B(p,2r)\cap\Omega}.
$$
On the other hand:
$$
\int_{\Omega}\abs{u}^2\geq \abs{B(p,r)\cap\Omega}.
$$
Hence its Rayleigh quotient is bounded above as follows:
$$
R(u)\leq\dfrac 1{r^2}\dfrac{\abs{B(p,2r)\cap\Omega}}{\abs{B(p,r)\cap\Omega}}\leq 
\dfrac{34}{r^2}.
$$
Recalling the definition of $r$, we conclude:
$$
R(u)\leq 544\pi\frac{n}{\vert \Omega \vert}
$$
as asserted.

\section{Proof of Theorem \ref{epsilonnet}}\label{partition}
 
The strategy of the proof is to partition the given punctured domain in a family of convex domains with only one puncture, and then to apply the lower bound \eqref{lbannulus} to each piece of the partition. 

\smallskip

First, we say that the family of open sets
$\{\Omega_1,\dots,\Omega_n\}$ with piecewise-smooth boundary is a {\it partition} of the open set $\Omega$ if 
$\bar\Omega=\cup_{j=1}^n\bar\Omega_j$; the partition is {\it disjoint} if moreover $\Omega_j\cap\Omega_k$ is empty whenever $j\ne k$. It is a simple consequence of the min-max principle that the first eigenvalue of $\Omega$ is controlled from below by the smallest first eigenvalue of the members of a disjoint partition, that is:
\begin{equation}\label{dp}
\lambda_1(\Omega,A)\geq\min_{j=1,\dots,n}\lambda_1(\Omega_j,A),
\end{equation}
for any potential one-form $A$ (for the easy proof we refer to  Proposition 4 of \cite{CS2}). The second ingredient is the estimate \eqref{lbannulus} for an annulus $\Omega=F\setminus \bar G$ with $F$ and $G$ convex with piecewise-smooth boundary:
\begin{equation}\label{bis}
\lambda_1(\Omega,A)\geq\dfrac{4\pi^2}{\abs{\bd F}^2}\dfrac{\beta^2}{B^2}d(\Phi,{\bf Z})^2
\end{equation}
where $\Phi$ is the flux of $A$ across the inner boundary component $\bd G$ and $\beta$ (resp. $B$) is the minimal and maximal width of $\Omega$, respectively.

\smallskip

Let us then start from the partition. In fact, 
the properties of a maximal $\epsilon$-net allow to partition the given domain in "well-balanced" convex pieces, in the following sense. 

\begin{lemme}\label{epartition} Let $\Omega$ be a convex domain with smooth boundary and let $\mathcal P_{\epsilon}=\{p_1,\dots,p_n\}$ be a maximal $\epsilon$-net in $\Omega$. We assume that $\bd\Omega$ satisfies the $\delta$-interior ball condition with $\delta>\epsilon$. Then $\Omega$ admits a disjoint partition $\{\Omega_1,\dots,\Omega_n\}$ with the following properties:

\smallskip

a) Every $\Omega_j$ is convex and has piecewise-smooth boundary;

\smallskip

b) For each $j=1,\dots,n$ one has $B(p_j,\frac{\epsilon}2)\subseteq\Omega_j\subseteq B(p_j,2\epsilon)$.
\end{lemme}

We will prove the lemma below. 

\smallskip

To finish the proof of Theorem \ref{epsilonnet}, we first observe that 
$\{\Omega_1\setminus\{p_1\},\dots,\Omega_n\setminus\{p_n\}\}$ is a disjoint partition of the punctured domain $\Omega\setminus\mathcal P_{\epsilon}$ and, in view of \eqref{dp},  it is enough to bound from below the ground state energy of every piece of this partition. To that end, we apply \eqref{bis} to  $\Omega_j\setminus\{p_j\}$, more precisely, we take $F=\Omega_j$, $G=B(p_j,\eta)$ and let $\eta\to 0$. Taking into account Lemma \ref{epartition}, we have, after taking the limit as $\eta\to 0$:
$$
\beta\geq \frac{\epsilon}2, \quad B\leq 2\epsilon,\quad \text{hence}\quad \frac{\beta}{B}\geq \frac 14
$$
because $\beta$ and $B$ tend, respectively, to the minimum and maximum distance of $p_j$ to $\bd\Omega_j$. 
Moreover, as $\Omega_j$ is convex, contained in $B(p_j,2\epsilon)$, we have by the monotonicity of the perimeter: $\abs{\bd\Omega_j}\leq 4\pi\epsilon$. The conclusion is that, for all $j=1,\dots,n$:
$$
\lambda_1(\Omega_j\setminus\{p_j\},A)\geq\dfrac{1}{64\epsilon^2}d(\Phi,{\bf Z})^2.
$$
As this holds for any member of the partition, it holds a fortiori for $\Omega\setminus\mathcal P_{\epsilon}$, which proves the first part of the theorem. 

Finally, it is readily seen that the number of points in a maximal $\epsilon$-net grows proportionally to $\epsilon^{-2}$. Precisely, one first observes that
$
\cup_{j=1}^{ n}B(p_j,\frac{\epsilon}2)\subseteq\Omega;
$
since the union on the left is disjoint (by maximality of the net) we obtain: $n\cdot\frac{\pi\epsilon^2}{4}\leq\abs{\Omega},$ that is
$$
\dfrac{1}{\epsilon^2}\geq\dfrac{\pi n}{4\abs{\Omega}},
$$
which proves \eqref{wein}.

\subsection{Proof of Lemma \ref{epartition}}

We will use the following property of maximal $\epsilon$-nets:

\medskip

{\bf Property P.} {\it If $x\in\Omega$ is such that $d(x,p_j)>\epsilon$ and $d(x,\bd\Omega)\geq\epsilon$, then there exists $p_k\ne p_j$ such that $d(x,p_k)\leq\epsilon$.}

\smallskip

 For each $j=1,\dots,n$ we consider the non-empty open set:
$$
\Omega_j=\{x\in \Omega: d(x,p_j)<d(x,p_k)\quad\text{for all}\quad k\ne j\}.
$$
It is clear that
$
\bar\Omega=\cup_{j=1}^n\bar\Omega_j.
$
If, for indices $j\ne k$ we consider the open half-space
\begin{equation}\label{hjk}
H_{jk}=\{x\in\real 2: d(x,p_j)<d(x,p_k)\}
\end{equation}
we see that we can write
$$
\Omega_j=\cap_{k\ne j}(H_{jk}\cap \Omega)
$$
which makes it clear that $\Omega_j$ is convex. As the boundary of $\Omega_j$ is either part of $\bd\Omega$, or is part of  $\bd H_{jk}$, which is a straight line, 
we see that $\bd\Omega_j$ is piecewise-smooth. This proves a). We now prove 
 the first inclusion in b). Assume $d(x,p_j)<\frac{\epsilon}{2}$: it is enough to show that $x\in\bar\Omega_j$. In fact, if $x\notin\bar\Omega_j$ there exists $k\ne j$ such that $x\in\bar\Omega_k$ and, by definition,  $d(x,p_k)\leq d(x,p_j)<\frac{\epsilon}{2}$. This means that
$$
d(x,p_j)< \frac{\epsilon}{2}\quad\text{and}\quad d(x,p_k)< \frac{\epsilon}{2}
$$
which by the triangle inequality gives $d(p_j,p_k)< \epsilon$, which  is a contradiction. Hence $x\in\bar\Omega_j$.

\medskip

We now prove the second inclusion in b). Let $x\in\Omega_j$. It is enough to show that $d(x,p_j)<2\epsilon$  in any of the following two cases: 

\smallskip

{\bf Case I}: $d(x,\bd\Omega)\geq\epsilon$,

\smallskip

{\bf Case II}: $d(x,\bd\Omega)<\epsilon$.

\smallskip

In Case I, assume  that $d(x,p_j)\geq 2\epsilon$, so that, in particular, $d(x,p_j)>\epsilon$. By Property P above there exists $p_k\ne p_j$ such that $d(x,p_k)\leq\epsilon$. As $x\in\Omega_j$ we have, by definition, $d(x,p_j)<d(x,p_k)$ hence a fortiori $d(x,p_j)<\epsilon$ which is a contradiction.

\medskip

Now assume we are in Case II. Let $\bar x\in\bd\Omega$ be the foot of the unique geodesic segment $\gamma$ which minimizes distance to the boundary, and let $p$ be a point of $\gamma$ at distance $\epsilon$ to $\bar x$. Since $\bd\Omega$ has the $\delta$-interior ball condition, and since $\epsilon<\delta$, it is clear that
$$
d(p,\bd\Omega)=\epsilon \quad\text{and}\quad d(p,x)\leq\epsilon.
$$
Since in particular $d(p,\bd\Omega)\geq\epsilon$, there exists $p_k\in\mathcal P_{\epsilon}$ such that $d(p,p_k)\leq\epsilon$, by the $\epsilon$-maximality of the net. By the triangle inequality:
$$
d(x,p_k)\leq d(x,p)+d(p,p_k)\leq 2\epsilon.
$$
On the other hand $x\in\Omega_j$ and because of that one has
$d(x,p_j)<d(x,p_k)$. Hence $d(x,p_j)<2\epsilon$ and the proof is complete.

%%%%%%%%%
%%%%%%%%%

\section{Bound of $\lambda_1(\Omega,A)$ with respect to the invariant $h(\Omega)$}

 The goal of this section is to prove Theorem \ref{main2} and to construct Example \ref{example2}. This example show that, surprisingly, when the number of holes increase, the constant $h$ can decrease to $0$, and, at the same time, the ground state energy can increase to $\infty$.

\subsection{Proof of Theorem \ref{main2}.}
 
 \begin{proof} 
 First, we observe that, by hypothesis, we can cut $n$ segments $\Gamma_1,...,\Gamma_n$ in $\Omega$ so that the complement
$\Omega \setminus \{\Gamma:=\Gamma_1 \cup...\cup \Gamma_n\}$
is simply connected.

\smallskip
Let $\Gamma(\epsilon)$ be the $\epsilon$-neighborhood of $\Gamma$ and set $D=\Omega\setminus \Gamma(\epsilon)$; if $\epsilon$ is small enough $D$ is simply connected and we have, by Proposition 12 in \cite{CS1}:  
$$
\lambda_1(\Omega,A) \le \nu_1(D)
$$
where $\nu_1(D)$ denotes the first eigenvalue of a mixed problem on $D$, where we take the Dirichlet condition on $\bd D\cap \Omega$ and the Neumann condition on $\partial D\cap \partial \Omega$.

In order to control $\nu_1(D)$, we will construct a test function $u$ taking the value $0$ on $\bd D\cap \Omega$ and apply the min-max principle.

\medskip
\noindent
 We assume $\pi h= \pi\sum_{i=1}^n h_i \le   \frac{\vert \Omega\vert}{2} $ and $h_i\leq 1$.  

\medskip
We consider one of the segments $\Gamma_i$ of length $h_i$, and denote by $q_i$ the middle of $\Gamma_i$.   

\begin{figure}[H]
\begin{center}
 \includegraphics[width=8cm]{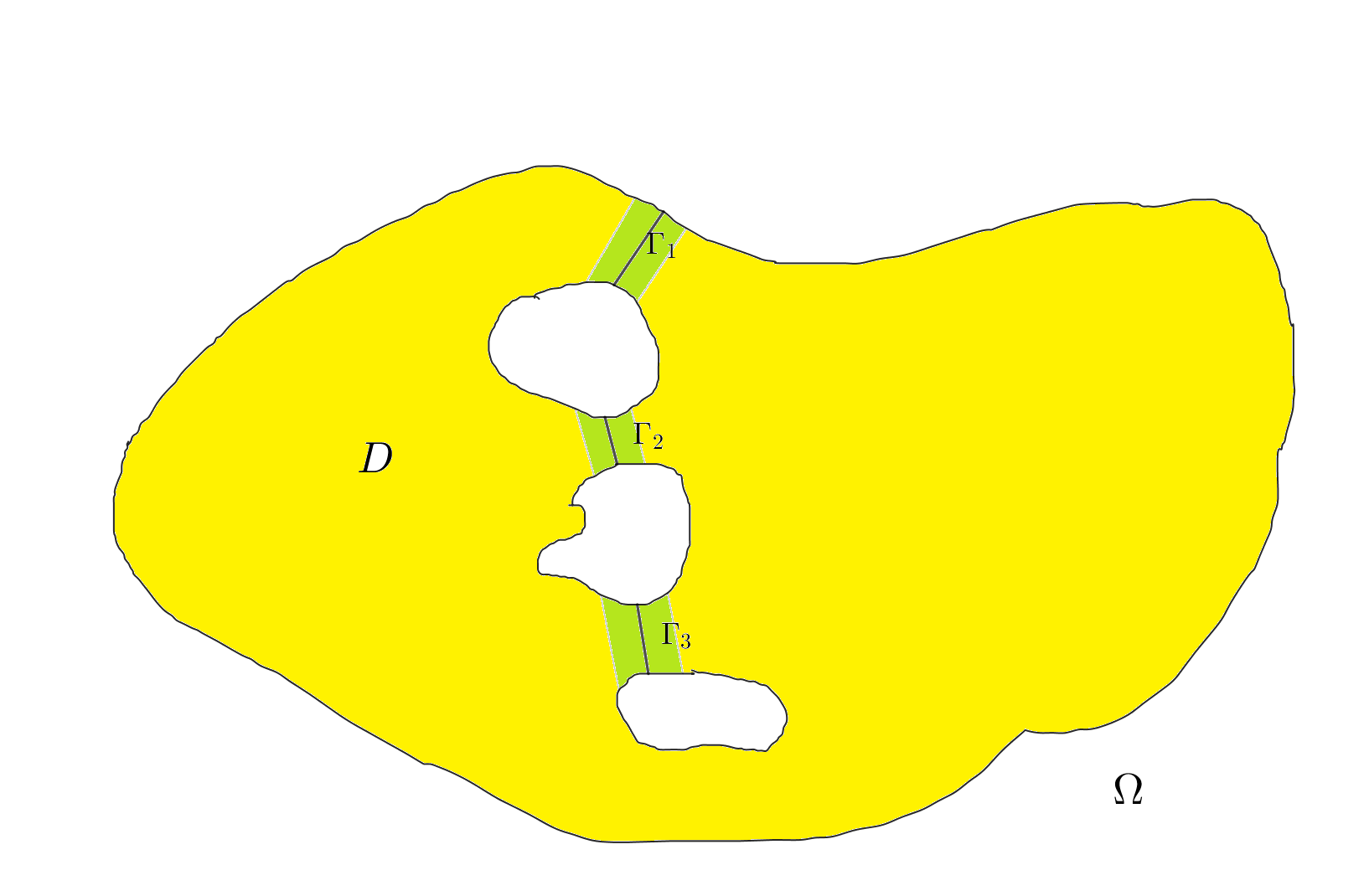}
\includegraphics[width=6cm]{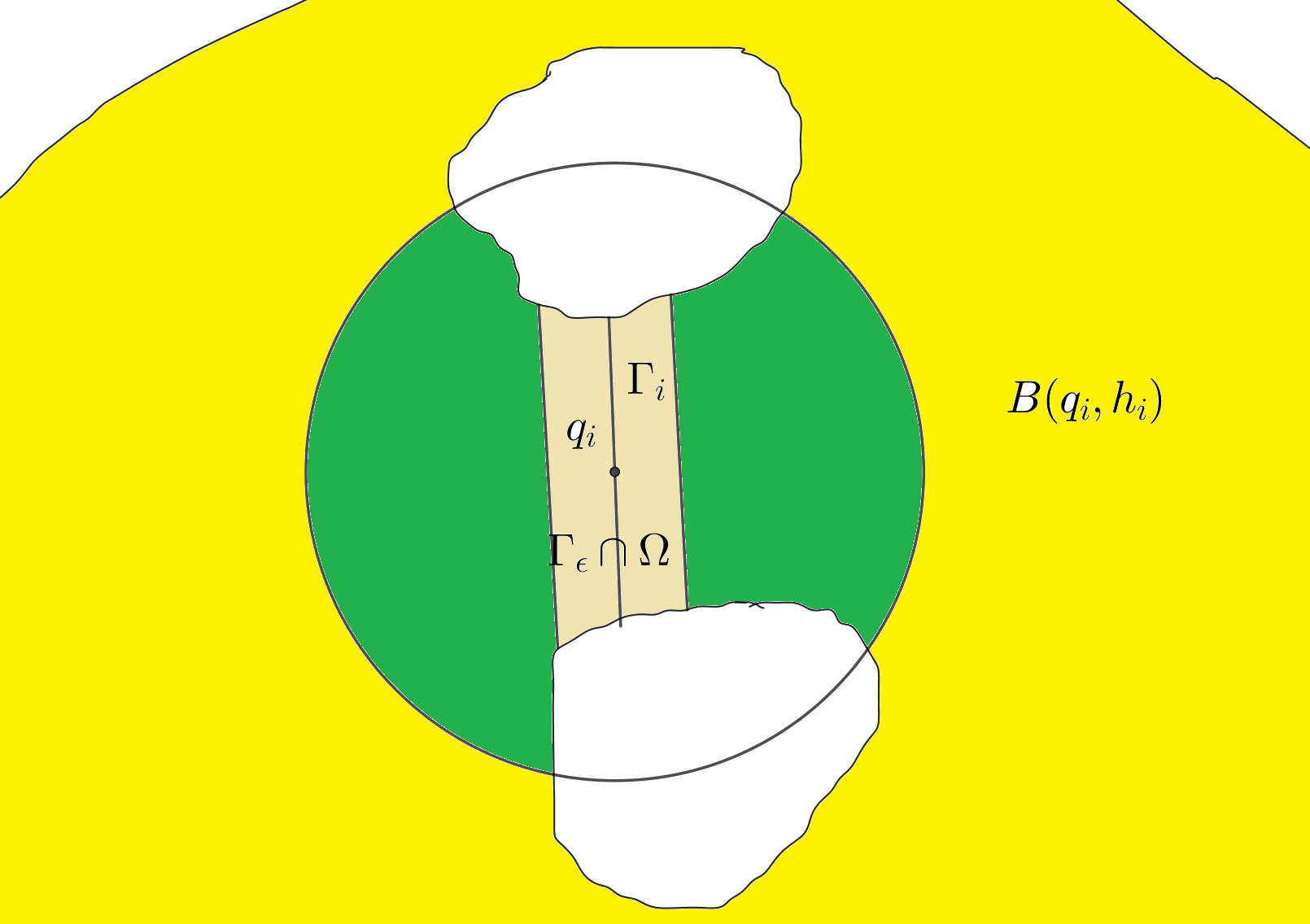}
\caption{On the left the domain $D= \Omega \setminus \Gamma_{\epsilon}$. On the right, the ball $B(q_i,h_i)$ and its intersection with $\Gamma_{\epsilon}$.}
\end{center}
 \end{figure}

\medskip
Observe that for $\epsilon $ small enough, $\Gamma(\epsilon)\subset B(q_1,h_1)\cup...\cup B(q_n,h_n)$. We will construct a test-function $u$ taking the value $0$ on $B(q_1,h_1)\cup...\cup B(q_n,h_n)$, so that it takes the value $0$ on $ \partial D \cap \Omega$.

\medskip
We introduce the radial function $u_i$ on $D$ defined by
$$
u_i(x)=\threesystem
{0\quad\text{if}\quad d(x,q_i)\leq h_i,}
{\frac{-2}{\ln h_i} (\ln d(x,q_i)-\ln h_i)\quad\text{if}\quad  h_i\le d(x,q_i)\le \sqrt {h_i},}
{1\quad\text{if}\quad d(x,q_i)\geq\sqrt{h_i}.}
$$

Our test function $u$ on $D$ will be the product:
$$
u(x)=u_1(x)\cdots u_n(x).
$$
Taking into account that $\vert u_i(x) \vert \le 1$, we have
$$
\vert \nabla u \vert  \le \vert \nabla u_1\vert +...+\vert \nabla u_n\vert,
$$
and therefore
$$
\vert \nabla u \vert^2 \le n (\vert \nabla u_1\vert^2 +...+\vert \nabla u_n\vert^2).
$$

So it suffices to bound from above the contribution of each $\int_{D}\vert \nabla u_i\vert^2$ in order to control $\int_D\vert \nabla u\vert^2$.
In polar coordinates centered at $q_i$, one has:
$$
\abs{\nabla u_i}^2=\dfrac{4}{\ln^2 h_i}\cdot\dfrac{1}{r^2}
$$
on the subset where $h_i\le r\le \sqrt {h_i}$, and zero everywhere else. Then:
$$
\int_{D}\vert \nabla u_i \vert^2\le 2\pi \frac{4}{\ln^2 h_i} \int_{h_i}^{\sqrt {h_i}}\frac{dr}{r}
=\frac{8 \pi}{\ln^2 h_i} (\ln \sqrt {h_i}-\ln h_i)= \frac{-4 \pi}{\ln h_i}
$$
hence summing over $i$ (recall that $h_i\leq 1$): 
$$
\int_D\vert \nabla u\vert^2 \le 4\pi n \sum_{i=1}^n \frac{-1}{\ln h_i/2}=
4\pi n \sum_{i=1}^n \frac{1}{\abs{\ln h_i/2}}
$$
Taking into account that the area of a ball of radius $\sqrt{h_i}$ is $\pi h_i$ and that $u=1$ outside these balls, the $L^2-$norm of the function $u$ is at least
$$
\int_{\Omega}u^2 \ge \vert \Omega\vert -\pi\sum_{i=1}^n h_i \ge \frac{\vert \Omega \vert}{2}
$$
because we assume $\sum_ih_i=h\leq\frac{\abs{\Omega}}{2\pi}$. So, by the min-max principle we deduce
$$
\lambda_1(\Omega,A)\le \frac{8\pi n}{\vert \Omega \vert}\sum_{i=1}^n \frac{1}{\abs{\ln h_i/2}}$$
as asserted.
\end{proof}

{\bf Proof of Corollary \ref{corh}.}
First observe that 
$$
\dfrac{1}{\ln{\frac 2{h_j}}}\leq -\dfrac{1}{\ln h_j}=\dfrac{1}{\abs{\ln h_j}}
$$
Then, since the function
$
\phi(x)=-\frac{1}{\ln x}
$
is concave ($\phi''(x)\leq 0$) on the interval $(0,e^{-2})$, we have,  by Jensen inequality:
$$
\sum_{j=1}^n\phi(h_j)\leq n\phi\Big(\dfrac{1}{n}\sum_{j=1}^nh_j\Big).
$$
Translated to our situation, we see that if every $h_j<e^{-2}$, then:
$$
\sum_{j=1}^n\dfrac{1}{\abs{\ln{\frac 2{h_j}}}}\leq \dfrac{n}{\abs{\ln(\frac{h(\Omega)}{n})}}
$$
and the upper bound \eqref{in4} reads:
$$
\lambda_1(\Omega,A)\le \frac{8\pi n^2}{\abs{\Omega}}\dfrac{1}{\abs{\ln(\frac{h(\Omega)}{n})}}.
$$
where $h(\Omega)$ is the Cheeger constant introduced in Definition \ref{Cheeger}.

\subsection{Construction of an example} \label{example}

We will construct a family of domains $\Omega_k$ with $n(\Omega_k)=k^2$ holes; each $\Omega_k$ is obtained as a union of $k^2$ identical fundamental pieces $C_k$. Each fundamental piece $C_k$ is a doubly-convex domain, so that we will be able to use the inequality of Theorem 1 in \cite{CS2} to bound from below its first eigenvalue.

\smallskip

\noindent\textbf{Step 1: the definition of the fundamental piece $C_k$.}

\smallskip
\noindent
The domain $C_k$ is determined by the exterior boundary curve, a square of sidelength $\frac 4k$ based on the vertices 
$$
A=\Big(-\frac{2}{k},0\Big), \quad B=\Big(\frac{2}{k},0\Big),\quad C=\Big(\frac{2}{k},\frac{4}{k}\Big), \quad D= \Big(-\frac{2}{k},\frac{4}{k}\Big), \quad 
$$
and the inner boundary curve, a rectangle based on the vertices
$$
A'=\Big(-\frac{1}{k},\frac{1}{k^{5/2}}\Big), \quad B'=\Big(\frac{1}{k},\frac{1}{k^{5/2}}\Big),\quad C'=\Big(\frac{1}{k},\frac{4}{k}-\frac{1}{k^{5/2}}\Big), \quad D'=\Big(-\frac{1}{k},\frac{4}{k}-\frac{1}{k^{5/2}}\Big).
$$

\smallskip

We refer to the picture below.

\begin{figure}[H]
\begin{center}
 \includegraphics[width=8cm]{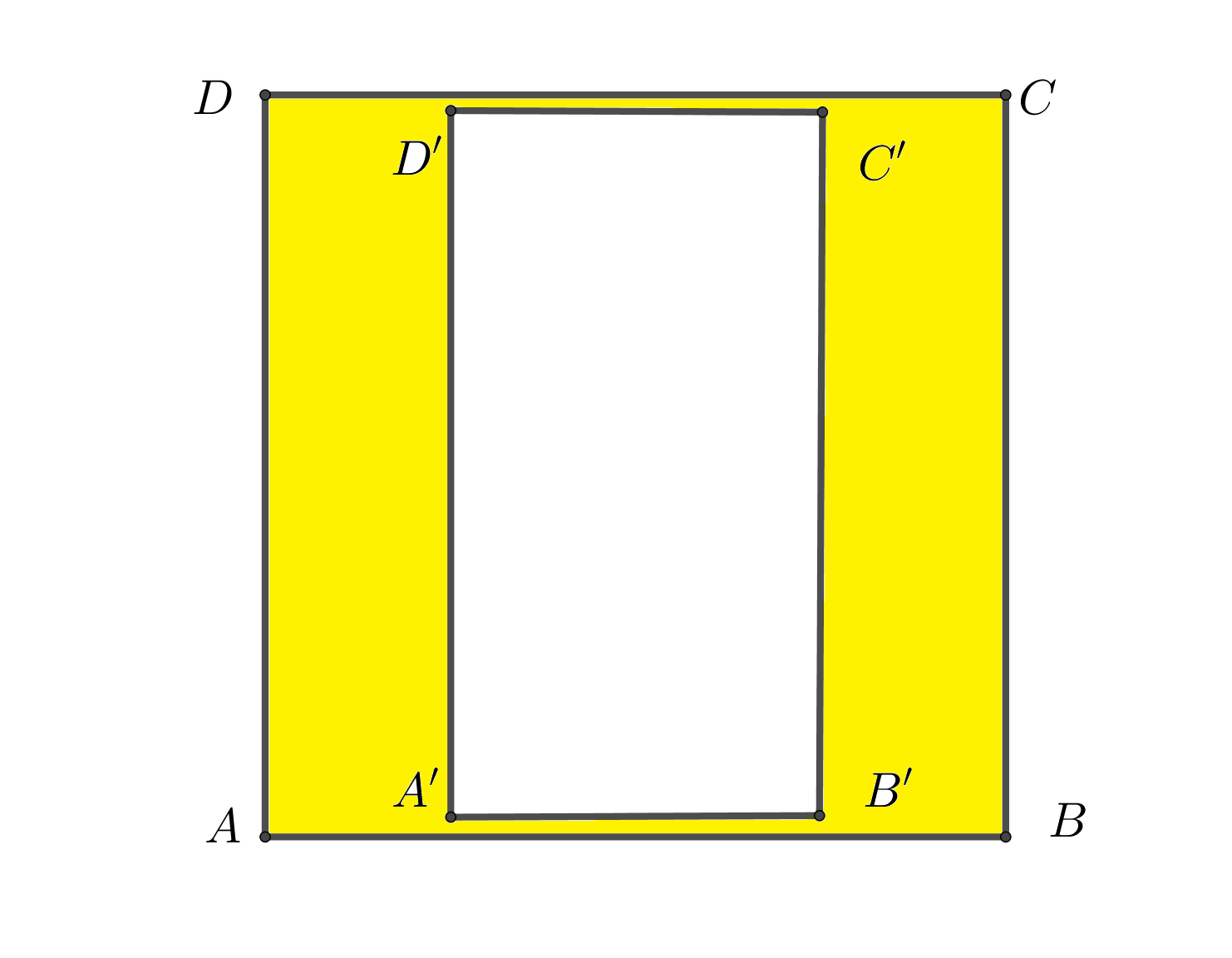}
\includegraphics[width=6cm]{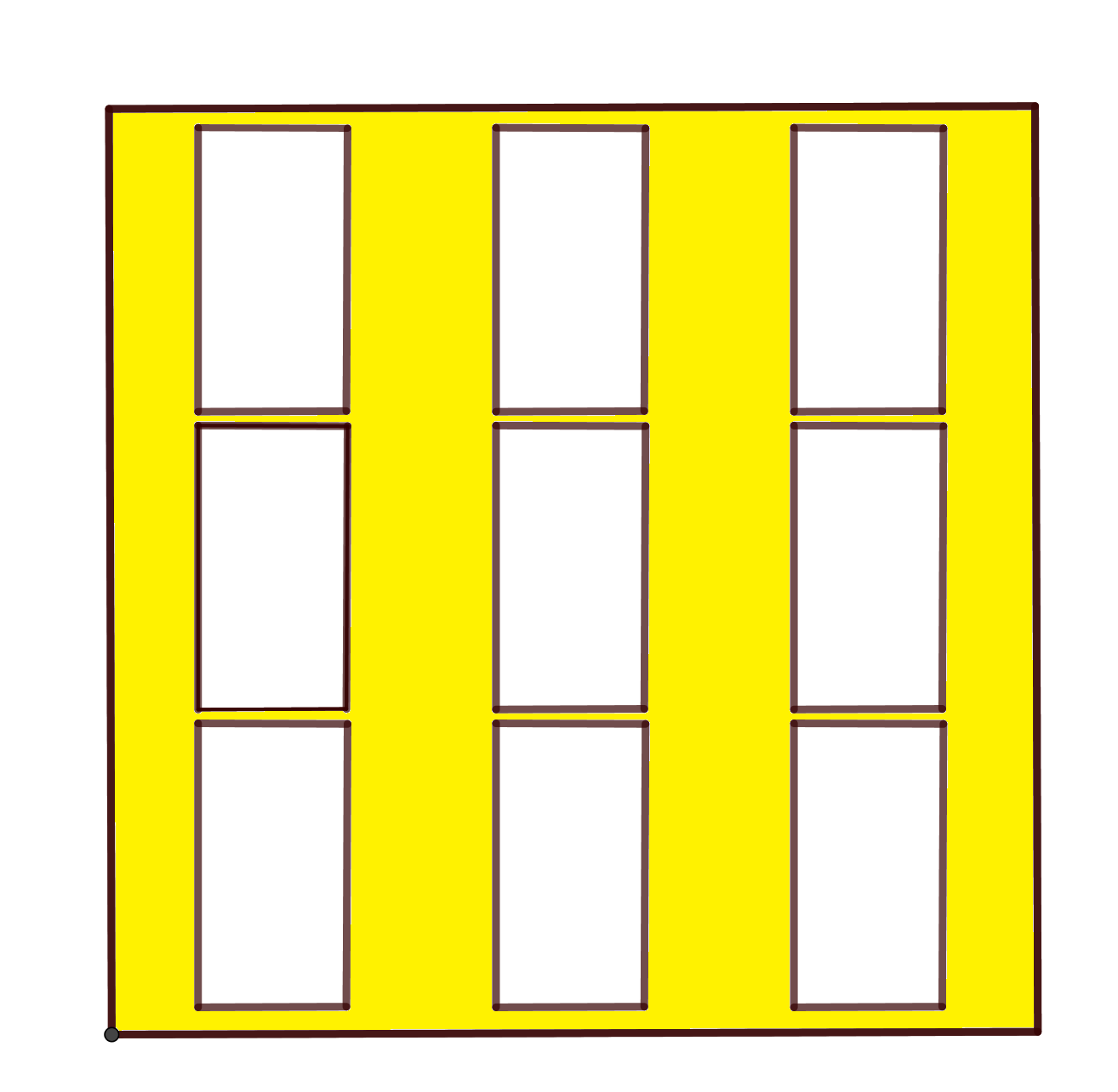}
\caption{On the left, the fundamental piece $C_k$ constructed with the rectangles $ABCD$ and $A'B'C'D'$. On the right, the domain $\Omega_k$ for $k=3$ obtained by assembling 9 fundamental pieces together.}
 \end{center}
\end{figure}

\noindent
 
We have:

\medskip
- the area of $C_k$ is $\vert C_k\vert= \frac{8}{k^2}+\frac{4}{k^{7/2}}$;   

\smallskip
- the area of the domain $F_k$ bounded by the outer curve is $\frac{16}{k^2}$;

\smallskip
- the minimal width $\beta_k$ between the two boundaries is $\beta_k= \frac{1}{k^{5/2}}$;

\smallskip
- the maximal width $B_k$ between the two boundaries is $B_k= \frac{1}{k}\sqrt{1+\frac{1}{k^3}}$;

\smallskip
- the diameter of $F_k$ is $D_k= \frac{4\sqrt 2}{k}$;

\smallskip
- the length $L_k$ of the exterior boundary is $L_k=\frac{16}{k}$.

\medskip
We consider a potential $A$ with flux $\Phi$. We can apply again Theorem 1 in \cite{CS2} (see \eqref{sharp}) and we obtain:

\begin{equation} \label{minoration}
\lambda_1(C_k,A)\geq c\, d(\Phi,{\bf Z})^2\sqrt k, \quad\text{with}\quad c=\frac{\pi^2}{2^{15}\sqrt{1+\frac{1}{k^3}}}\geq \dfrac{\pi^2}{2^{15}\sqrt 2},
\end{equation}
which grows like $\sqrt k$.

\medskip
\noindent
\textbf{Step 2: the definition of $\Omega_k$ and its first eigenvalue.} The domain $\Omega_k$ is a square of size $4$ filled with $k^2$ copies of $C_k$ as in the picture. That is:
$$
\Omega_k=\cup_{j=1}^{k^2}C_{kj}
$$
where $C_{kj}$ is congruent to $C_k$ for all $j$. 
$\Omega_k$ has $k^2$ holes and $\vert \Omega \vert= k^2\vert C_k\vert= 8+\frac{4}{k^{3/2}} \ge 1$.

 \medskip
We consider a potential $A$ having fixed flux $\Phi$ around each hole.
First, it is easy to show that
$$
\lambda_1(\Omega_k,A)\ge \lambda_1(C_k,A)\geq cd(\Phi,{\bf Z})^2 k^{\frac{1}{2}}
$$
where $c$ is the constant in Formula (\ref{minoration}). In fact, let $u$ be a first eigenfunction of $\Omega_k$. By restricting it to each piece $C_{kj}$ we have, by the min-max principle:
$$
\int_{C_{kj}}\abs{\nabla^Au}^2\geq  \lambda_1(C_{kj},A)\int_{C_{kj}}\abs{u}^2= \lambda_1(C_k,A)\int_{C_{kj}}\abs{u}^2.
$$
We sum over $j=1,\dots,k^2$ and obtain
$$
\int_{\Omega_k}\abs{\nabla^Au}^2\geq\lambda_1(C_{k},A)\int_{\Omega_k}\abs{u}^2
$$
which immediately gives $\lambda_1(\Omega_k,A)\geq \lambda_1(C_k,A)$.

 \medskip
\noindent
\textbf{Step 3: calculation of $h(\Omega_k)$.}
Now, let us see the total length $h$ of the segments we have to cut in order to make $\Omega_k$ simply connected.
As the holes are at distance $\frac{2}{k^{5/2}}$, we need to cut $(k-1)k$ segments of length $ \frac{2}{k^{5/2}}$ and $k$ segments of length $\frac{1}{k^{5/2}}$. 

\smallskip

\noindent The total length we need to cut is $h(\Omega_k)= \frac{2k-1}{k^{3/2}}$. Summarizing we have, as $k\to\infty$:
$$
\abs{\Omega_k}\ge 1, \quad h(\Omega_k)\sim \frac{2}{\sqrt k}, \quad\lambda_1(\Omega_k,A)\geq c \, d(\Phi,{\bf Z})^2\sqrt k,
$$
with $ c\geq \frac{\pi^2}{2^{15}\sqrt 2}$, which in particular shows that $\abs{\Omega_k}$ is bounded from below, $h(\Omega_k)$ tends to zero and $\lambda_1(\Omega_k,A)$ tends to infinity, as requested.

\addcontentsline{toc}{chapter}{Bibliography}
%\listoffigures
%\nocite{*} %Permet d'afficher toute la bibliographie sans forcément avoir fait référence dans le document
\bibliographystyle{plain}
\bibliography{biblioCSdomains1}

 \bigskip
\normalsize 
\noindent Bruno Colbois \\
Universit\'e de Neuch\^atel, Institut de Math\'ematiques \\
Rue Emile Argand 11\\
 CH-2000, Neuch\^atel, Suisse

\noindent bruno.colbois@unine.ch

\medskip

\normalsize
\noindent
Alessandro Savo \\
 Dipartimento SBAI, Sezione di Matematica,
Sapienza Universit\`a di Roma\\
Via Antonio Scarpa 16\\
00161 Roma, Italy

\noindent alessandro.savo@uniroma1.it

\end{document}